\input amstex
\documentstyle{amsppt}

\document
\vsize=6 in \hsize=5.5 in \magnification=1200

\topmatter
\abstract {In this note we analyse the Exceptional Component of the space of integrable forms of degree two, introduced by Cerveau-Lins Neto, in terms of the geometry of Veronese curves and classical invariant theory.} \endabstract
\title  A note on the $\jmath$ invariant and foliations   \endtitle
\author  Omegar  Calvo-Andrade and  Fernando Cukierman
\endauthor
 \thanks{We thank Jorge Vit\'orio Pereira for some useful communications.} \endthanks
\address {CIMAT, Ap. postal 402, Guanajuato 36000, Mexico}  \endaddress
\email {omegar\@cimat.mx} \endemail
\address {Dto. Matematica, FCEN-UBA,
Ciudad Universitaria, (1428) Buenos Aires, Argentina}
\endaddress
\email {fcukier\@dm.uba.ar} \endemail
\date{ }\enddate
\endtopmatter

\comment
\NoRunningHeads
\NoPageNumbers
\leftheadtext{  }
\rightheadtext{  }

\flushpar
{\bf Contents} \newline
\S 1   \newline
\S 2  \newline
\S 3 \newline
\S 4 \newline
\S 5  \newline
\S 6  \newline
\S 7
\endcomment

\flushpar
{\bf 1. Introduction}
\newline
\newline
(1.1) Let $r$ and $d$ be natural numbers. Consider a differential 1-form in $\Bbb C^{r+1}$
$$\omega = \sum_{i=0}^r a_i  dx_i$$
where the $a_i$ are homogeneous polynomials of degree $d + 1$ in
variables $x_0, \dots, x_r$, with complex coefficients. Let us
assume that
$$ \sum_{i=0}^r a_i  x_i = 0$$
so that $\omega$ descends to the complex projective space $\Bbb P^r$
and defines a global section of the twisted sheaf of 1-forms
$\Omega^1_{\Bbb P^r}(d+2)$.
\newline
\newline
(1.2) For a $K$-vector space $V$ we denote $\Bbb P V = V - \{0\} /K^*$ the
projective space of one-dimensional linear subspaces of $V$
and $\pi: V - \{0\} \to \Bbb P V$  the quotient map.
\newline
\newline
Consider the  projective space $\Bbb P H^0(\Bbb P^r, \Omega^1_{\Bbb
P^r}(d+2))$ and the subset
$$\Cal F(r, d) = \pi (\{ \omega \in
H^0(\Bbb P^r, \Omega^1_{\Bbb P^r}(d+2)) -\{0\}
 / \ \omega \wedge d\omega = 0\})$$
parametrizing 1-forms $\omega$ that satisfy the Frobenius integrability condition.
\newline
\newline
(1.3) It is clear that $\Cal F(r, d)$ is an algebraic subset defined
by quadratic equations. It is the space of degree $d$ foliations of
codimension one on $\Bbb P^r$. One important problem of the area is
to determine the irreducible components of $\Cal F(r, d)$.
\newline
\newline
To fix notation, let us recall (see \cite{CL}) the following families of irreducible components:
\newline
\newline
a) The rational components $\Cal R(d_1, d_2) \subset \Cal F(r, d)$
consisting of the 1-forms of type
$$\omega = p_1 \ F_2 \ dF_1 - p_2 \ F_1 \ dF_2$$
where $d + 2 = d_1 + d_2$ is a partition with $d_1, d_2$ natural
numbers, $p_1, p_2$ are the unique coprime natural numbers such that
$p_1 d_1 = p_2 d_2$ and $F_1, F_2$ are homogeneous polynomials of
respective degrees $d_1, d_2$.
\newline
\newline
b) The logarithmic components $\Cal L(d_1, \dots, d_s) \subset \Cal F(r, d)$
consisting of the 1-forms of type
$$\omega = \left(\prod_{j=1}^s  F_j\right) \sum_{i=1}^s \lambda_i \ dF_i /F_i =
 \sum_{i=1}^s \left(\prod_{j \ne i}  F_j \right) \ \lambda_i \ dF_i $$
where $s \ge 3$, the $d_i$ are integers such that $d + 2 =
\sum_{i=1}^s d_i$, the $F_i $ are homogeneous polynomials of degree
$d_i$ and $\lambda_i$ are complex numbers, not all zero, such that
$\sum_{i=1}^s d_i \lambda_i = 0$.
\newline
\newline
c) The linear pull-back components $PBL(d) \subset \Cal F(r, d)$
consisting of the 1-forms of type
$$\omega = \pi^* \eta$$
where $\pi:\Bbb C^{r+1} \to \Bbb C^3$ is a non-degenerate linear map
and $\eta$ is a global section of $\Omega^1_{\Bbb P^2}(d+2)$.
\newline
\newline
(1.4) The problem of determining the irreducible components of $\Cal
F(r, d)$ was solved by D. Cerveau and A. Lins Neto in \cite{CL} for
$d=2$. These authors defined an irreducible component $\Cal E
\subset \Cal F(3, 2)$, called "the exceptional component". The
leaves of a typical foliation in $\Cal E$ are the orbits of a linear
action in $\Bbb P^3$ of the affine group in one variable. The main
theorem of \cite{CL} states that the irreducible components of $\Cal
F(r, 2)$ are the corresponding rational, logarithmic and linear
pull-backs components and the component $\Cal E_r$ obtained by
linear pull-backs $\Bbb P^r \to \Bbb P^3$ of $\Cal E = \Cal E_3$.
\newline
\newline
(1.5) The purpose of this note is to give another description of
$\Cal E$, emphasizing the role of the $\jmath$ invariant of a binary
quartic; see also \cite{CD}, example (2.4.8), page 36. We shall
consider the  codimension one foliation on $\Bbb P^4$ induced by the
natural action of $PGL(2, \Bbb C)$ and will obtain $\Cal E$ by
restricting to a suitable hyperplane $\Bbb P^3 \subset \Bbb P^4$,
namely, an osculating hyperplane of the Veronese curve.
\newline
\newline
(1.6) We introduce here some notation that will be useful later. Let
$\omega \in H^0(\Bbb P^r, \Omega^1_{\Bbb P^r}(d+2))$ as above.
Denote
$$S(\omega)$$ the variety of zeros of $\omega$ and $S_k(\omega)$ the
union of the irreducible components of $S(\omega)$ of dimension $k$.
If $S_{r-1}(\omega)$ is non-empty (i.e. if $\omega$ vanishes in
codimension one) then there exists a homogeneous polynomial $F$ of
maximal degree  $0 < e < d $ that divides  $\omega$. We denote
$$\bar \omega = \omega / F \in H^0(\Bbb P^r, \Omega^1_{\Bbb P^r}(d - e +2))$$
It is clear that $\bar \omega$ is well defined up to multiplicative constant and it does not vanish in codimension one.
\newline
\newline
\newline
{\bf 2.} Let us recall some known facts about the Veronese curves in $\Bbb P^4$,
following \cite{D} and \cite{H}.
\newline
\newline
(2.1) Our model for $\Bbb P^4$ will be the projective space of binary quartics.
\newline
\newline
Let $V$ be a two dimensional vector space over $\Bbb C$ and
for each natural number $r$ denote
$$P(r) = \Bbb P \text{Sym}^r(V)$$
the projective space associated to the $r+1$-dimensional vector space
$\text{Sym}^r(V)$.
\newline
\newline
The general linear group  $G = GL(V)$ acts on $V$ and hence naturally acts on $\text{Sym}^r(V)$ and on $P(r)$.
\newline
\newline
(2.2) Consider the Veronese map
$$\nu_r: P(1) \to  P(r)$$
obtained by sending $v \in V$ to $v^r \in \text{Sym}^r(V)$. It is clear that
$\nu_r$ is $G$-equivariant. The image $X_r = \nu_r(P(1))$ is called $r$-th Veronese curve. Notice that $G$ acts linearly on $P(r)$ and preserves $X_r$.
Hence $G$ also preserves the tangential and secant varieties of $X_r$.
\newline
\newline
To write this out in coordinates, let $t_0, t_1$ denote a basis of $V$.
Then $\{t_0^{r-i} t_1^i, \ i = 0,  \dots, r\}$ is a basis of $\text{Sym}^r(V)$
and a typical element of $\text{Sym}^r(V)$ is a binary form
$$F = \sum_{1=0}^r  a_i  t_0^{r-i}  t_1^i$$
Thus $\nu_r$ is defined by
$$\nu_r(a_0 t_0 + a_1 t_1) = (a_0 t_0 + a_1 t_1)^r =
\sum_{1=0}^r  \binom r i   a_0^{r-i}  a_1^i  t_0^{r-i}  t_1^i$$
By assigning to a homogeneous polynomial $F \in \text{Sym}^r(V)$
its roots counted with multiplicity, we may conveniently think of $P(r)$
as the set of effective divisors of degree $r$ in $P(1)$.
\newline
\newline
In these terms, a typical point in $X_r$ is a divisor of the form $r. p$
for some $p \in P(1)$. The tangent line to $X_r$ at this point is the set
of divisors of the form $(r-1) p + q$ for $q \in P(1)$. More
generally, the osculating $k$-plane to $X_r$ at $r. p$ is the set of divisors
of the form $(r-k)p + A$ where $A$ is any effective divisor of degree $k$
in $P(1)$.
\newline
\newline
(2.3) Let $r=4$. The orbits of $G$ in $P(4)$ are:
\newline
\newline
a) $X_4 = \{4p/ \ p \in P(1)\}$ = binary forms with a four-fold root
= Veronese curve. It is the unique closed orbit.
\newline
\newline
b) $T = \{3p + q/ \  p, q \in P(1), p \ne q \}$ = binary forms
with a triple root. The closure $\bar T$ equals the tangential surface
(union of all tangent lines) of $X_4$.
\newline
\newline
c) $N = \{2p + 2q/ \  p, q \in P(1), p \ne q \}$ =
binary forms with two double roots = set of points of intersection of pairs of
osculating 2-planes of $X_4$.
\newline
\newline
d) $\Delta = \{2p + q + r/  \  p, q, r  \text{ distinct points of }  P(1) \}$
= binary forms with one double root.
The closure $\bar \Delta$ is the discriminant
hypersurface, consisting of binary forms that have a multiple root.
It equals the union of all osculating 2-planes of $X_4$.
\newline
\newline
Notice that
$\bar T \cup \bar N \subset \bar \Delta$ and $\bar T \cap \bar N = X_4$.
\newline
\newline
e) Denote
$U = \{p + q + r + s/  \  p, q, r , s \text{ distinct points of }  P(1)  \}
= P(r) - \bar \Delta$
the open set of binary forms with simple roots. Then $U$ is a disjoint union of
infinitely many $G$-orbits. More precisely, one has the classical function
$$\jmath: U \to \Bbb C$$
such that the orbits are the sets $\jmath^{-1}(t)$ for $t \in \Bbb C$.
Thus, the foliation on $U$ with leaves the $G$-orbits is defined by
the differential 1-form $d\jmath$.
\newline
\newline
(2.4) As in \cite{D}, an explicit formula for $\jmath$ may be written in the form
$$\jmath = \frac {Q^3} {D}$$
where $Q, D$ are the homogeneous polinomials of respective degrees 2 and 6
given by :
$$
\align
Q &= a_0 a_4 - 4 a_1a_3 + 3 a_2^2 \\
C & = a_0 a_2 a_4 - a_0a_3^2 + 2a_1a_2 a_3 - a_1^2a_4 - a_2^3 \\
D &= Q^3 - 27 C^2
\endalign
$$
Here $D$ is the {\it discriminant} of a binary quartic,
so that $\bar \Delta $ is the set of zeros of $D$.
\newline
\newline
$Q$ is an equation for the unique $G$-invariant quadric containing $X_4$.
\newline
\newline
The cubic $C$ is called the {\it catalecticant} and is an
equation for the secant variety $\text {Sec }X_4$ of the Veronese curve.
\newline
\newline
It is also true that $\{Q, C\}$ generate the ring of invariants, but we will
not use this fact here.
\newline
\newline
We may consider $\jmath$ as a rational function $\jmath: P(4) \to \Bbb P^1$,
regular in the complement of the base locus $(C = Q = 0)$.
In other terms, $\jmath$ is the rational map to $\Bbb P^1$
defined by the pencil of sextic hypersurfaces in $P(4)$
spanned by $\{Q^3, D\}$ (or, equivalently, by $\{Q^3, C^2\})$.
Notice that all the sextics of this linear system are singular along the base locus.
\newline
\newline
(2.5) There are three fibers of $\jmath$ that deserve special attention:
$$
\align
& \jmath^{-1}(0) =  (Q=0)    \\
& \jmath^{-1}(\infty) =  (D=0) = \bar \Delta \\
& \jmath^{-1}(1728) =  (C=0) = \text {Sec }X_4
\endalign
$$
Taking account of multiplicities and writing $\jmath^*$ for
pull-back of divisors, we have:
$$
\align
& \jmath^*(0) =  3 (Q=0)   \\
& \jmath^*(\infty) =  (D=0)  \\
& \jmath^*(1728) =  2 (C=0)
\endalign
$$
The fiber at $\infty$ is reduced and irreducible, but it is singular in
codimension one. In fact,
$$\text{Sing }(\bar \Delta) = \bar T \cup \bar N$$
and, more precisely,
$\bar \Delta$ is cuspidal along $T$ and nodal along $N$.
\newline
\newline
Since each orbit is smooth and irreducible and
each fiber of $\jmath$ is a union of orbits, it
follows from the description of the orbits in (2.3)
that all other fibers of $\jmath$ are irreducible and smooth away from the base locus.
\newline
\newline
(2.6) Consider the codimension-one singular foliation $F$ in $P(4)$ with leaves the
fibers of $\jmath$, that is, $F$ is the singular foliation induced by the
natural action of $PGL(2, \Bbb C)$ on $P(4)$.
It will be convenient to consider the rational function
$$\jmath' = \frac {\jmath} {27(\jmath - 1)} = \frac {Q^3} {C^2}$$
Since $\jmath'$ and $\jmath$ differ only by an automorphism of $\Bbb P^1$,
they define the same foliation. Therefore,
the foliation  $F$ is defined by the differential form
$$\omega = 3 Q dC - 2 C dQ $$
and hence belongs to the irreducible component $\Cal R(2, 3) \subset
\Cal F(4, 3)$.
\newline
\newline
In order to describe the zeros of $\omega$,
we need another fact about the geometry of $X_4$.
\newline
\newline
(2.7) {\bf Proposition:} As in (2.3), let $\bar T$ denote the tangential surface
of the Veronese curve $X_4$. Then $\bar T = (C =0) \cap (Q = 0)$ and
the intersection is generically transverse. In particular, the base locus
of the pencil $\jmath$ is the tangential surface $\bar T$.
\newline
\newline
{\bf Proof: } We have $\bar T \subset (C =0)$ since in general the tangent variety is contained in the secant variety.
The inclusion $\bar T \subset (Q = 0)$ follows from direct calculation with
the formula for $Q$ in (2.4) or by  \cite {FH}, Ex. 11.32.
Then, $\bar T \subset (C =0) \cap (Q = 0)$.
On the other hand, it is shown in \cite {H}, p. 245, that the tangential surface
of the Veronese
curve $X_r \subset P(r)$ has degree $2r-2$. Therefore $\bar T$ has degree 6.
Since $(C =0) \cap (Q = 0)$ also has degree 6, the desired equality
and transversality hold.
\newline
\newline
(2.8) {\bf Proposition:}  $S(\omega) = \bar T \cup \bar N$.
In particular, all the irreducible components of $S(\omega)$ are of codimension two in $P(4)$.
\newline
\newline
{\bf Proof:} The zeros of $\omega$ consist of the base locus
$(C =0) \cap (Q = 0) = \bar T$ and of the singularities
of the fibers of $\jmath$. We know from (2.5) that the fiber
at $\infty$ is singular along $\bar N$ and the fibers
$\jmath^{-1}(t)$ are smooth away from the base locus for $t \notin \{0, 1728, \infty\}$;
this implies the result.
\newline
\newline
(2.9) Now we consider the restriction $F_H$ of $F$ to a hyperplane
$H \subset P(4)$. The singularities of $F_H$ are:
\newline
\newline
a) the intersections with $H$ of the singularities of $F$, and
\newline
\newline
b) the tangencies of $F$ and $H$ (that is, the loci of contact of the leaves of $F$ not transverse to $H$).
\newline
\newline
Denoting $\omega_H$ the 1-form in $H$ obtained by restriction of $\omega$, the foliation $F_H$ is defined by $\bar \omega_H$, with notation as in (1.6).
\newline
\newline
If $H$ is a general hyperplane then $\omega_H$ does not vanish in
codimension one. Hence $F_H$ is defined by $\omega_H$ and is a
rational foliation of type $\Cal R(2, 3)$ in $H \cong \Bbb P^3$. In
particular $F_H$ is a degree $3$ foliation.
\newline
\newline
(2.10) Now we analize $F_H$ when $H$ is an osculating hyperplane to the Veronese curve. Let $p \in X_4$ be a point and consider the osculating flag of $X_4$ at p:
$$\Bbb P^1_p = \{ 3p + q, \  q \in \Bbb P^1 \} \subset
\Bbb P^2_p = \{ 2p + q + r, \  q, r \in \Bbb P^1\} \subset
H = \Bbb P^3_p  = \{ p + q + r + s, \  q, r, s \in \Bbb P^1\}$$
Let us remark that the set
$$X_2 = \{2p + 2 q, \  q \in \Bbb P^1 \} \subset \Bbb P^2_p$$
is a copy of a Veronese curve of degree two in $\Bbb P^2$, and
$$X_3 = \{p + 3 q, \  q \in \Bbb P^1 \} \subset \Bbb P^3_p$$
is a copy of a Veronese curve of degree three in $\Bbb P^3$.
\newline
\newline
(2.11) {\bf Proposition:}  Let $\omega = 3 Q dC - 2 C dQ$ denote as above
the 1-form in $P(4)$ defining the foliation $F$ and $\omega_H$ its restriction to $H$. Then the zeros in codimension one of $\omega_H$ are
$$S_2(\omega_H) = \Bbb P^2_p \ \subset \ H$$
and the zeros of $\bar \omega_H$ (notation as in (1.6)) are
$$S(\bar \omega_H) = S_1(\bar \omega_H) =
\Bbb P^1_p \cup X_2 \cup X_3 \ \subset \  H.$$ In particular the
foliation induced by $\bar \omega_H$ has degree two.
\newline
\newline
{\bf Proof:} Following (2.9)a), let us determine
$S(\omega) \cap H = (\bar T \cup \bar N) \cap H$.
We find, set theoretically,
$$\bar T  \cap H = \{ 3r + s, \  r, s \in \Bbb P^1 \} \cap H =
\{3p + s, \  s \in \Bbb P^1 \} \cup \{ 3r + p, \  r \in \Bbb P^1 \} =
\Bbb P^1_p \cup X_3$$
$$\bar N \cap H = \{ 2r + 2s, \  r, s \in \Bbb P^1 \} \cap H =
\{ 2p + 2s, \  s \in \Bbb P^1 \} = X_2$$
\newline
Next, according to (2.9)b), let us look for tangencies.
We claim that the leaf $\jmath^{-1}(\infty) = \bar \Delta$ is tangent to $H$ along $\Bbb P^2_p$. In fact, the intersection $\bar \Delta  \cap H$ has two irreducible components, namely
$$\bar \Delta  \cap H = \{ 2q + r + s, \  q, r, s \in \Bbb P^1\} \cap H =
\{ 2p + r + s, \   r, s \in \Bbb P^1\} \cup
\{ 2q + r + p, \  q, r \in \Bbb P^1\} = \Bbb P^2_p \cup  \bar \Delta_H $$
where
$$\bar \Delta_H = \{ 2q + r + p, \  q, r \in \Bbb P^1\} \subset \Bbb P^3_p$$
is a copy of the discriminant hypersurface
$\{ 2q + r , \  q, r \in \Bbb P^1\} \subset \Bbb P(3)$ consisting of singular
cubic binary forms, and is hence an irreducible surface of degree four.
(In general, the discriminant for homogeneous polynomials of degree $d$ in
$n$ variables is an irreducible hypersurface of degree $n (d-1)^{n-1}$,
see \cite{GKZ}).
\newline
\newline
Since $\bar \Delta$ is the union of the osculating planes of $X_4$, it follows (see e. g. \cite{H}, Exercise (17.10)) that $\bar \Delta$ and $H$ are tangent to each other along $\Bbb P^2_p$, as claimed.
\newline
\newline
Notice that since $\bar \Delta$ is a sextic, we obtain the equality of divisors in $H$
$$\bar \Delta . H = 2  \Bbb P^2_p  + \bar \Delta_H$$
In particular $\bar \Delta$ and $H$ are transverse generically along
$\bar \Delta_H$ and therefore the only tangency contributed by the
leaf $\bar \Delta$ is the one along $\Bbb P^2_p$. To finish we only
need to see that there are no other leaves tangent to $H$. This
follows from the fact that the orbits of binary cubics under the
affine group are the same as before, Q.E.D.
\newline
\newline
\newline
{\bf 3.} We end this note with two  proofs, alternative to the one
given in \cite{CL}, of the fact that the closure of the orbit of
$\bar \omega_H$, which we now denote $\Cal E$, is an irreducible
component of $\Cal F(3, 2)$. We denote $F$ the foliation
defined by $\bar \omega_H$.
\newline
\newline
(3.1) Using the formulas in (2.4) we may write in appropriate coordinates
$$\bar \omega_H = x_3 [(2x_1^2 - 3 x_0 x_2) dx_0 + (3 x_2 x_3 - x_0 x_1) dx_1
+ (x_0^2 - 2 x_1 x_3) dx_2] - (x_0 x_1^2 - 2 x_0^2 x_2 + x_1 x_2
x_3) dx_3$$ as in \cite{CL}.

A straightforward computation shows that the singular set of
$d(\bar \omega_H)$ is  a point, namely,
the intersection point of the cubic, conic and line in
$S(\bar \omega_H)$.  It follows from Corollary 1 (section 5.2) or
Corollary 6.1 of \cite{CP} that every foliation $F'$ sufficiently
close to $F$ is induced by an action.
\newline
\newline
(3.2) It will be convenient to give explicit expressions for the
vector fields on $\Bbb P^3$ inducing $F$. If we write the
action of $\text{Aff}(\Bbb C)$ on $\Bbb C_3[t]$ as
$(at+b) \cdot p(t) = p(at +b)$ then the generators $x = (1+\epsilon)t$
and $y = t+\epsilon$ of $\text{Aff}(\Bbb C)$ act on basis
elements $t^i$ as follows:
$$
 x.t^i = ((1+\epsilon)t)^i  =  t^i  + i \epsilon t^{i}
\, \, \text{ and } \, \,
 y.t^i =  (t+\epsilon)^i    =  t^i + \epsilon i t^{i-1} \  \mod \epsilon^2
$$
It follows that the tangent sheaf of $F$ is generated by the vector fields
$$
X =   \sum_{i=0}^3 i z_i \frac{\partial}{\partial z_i} \, \, \text{
and }  \, \,  Y = \sum_{i=1}^{3} i z_{i-1} \frac{\partial}{\partial
z_i} .
$$
After a change of coordinates of the form $(z_0,\ldots, z_3) \mapsto
(\lambda_0 z_0,\ldots, \lambda_3 z_3)$ we can assume that
$$
X= \sum_{i=0}^3 i z_i \frac{\partial}{\partial z_i} \, \, \text{ and
} \, \, \,  Y= \sum_{i=1}^{3} z_{i-1} \frac{\partial}{\partial z_i}
.
$$
Notice that $[X,Y]=-Y$.
\newline
\newline
(3.3) Since the affine Lie algebra is rigid, the foliation $F'$
is induced by a $1$-form $\omega'$ of the form $i_{X'}i_{Y'}i_R
\Omega$ where $X'$, $Y'$ are close to $X$, $Y$ respectively and
$[X',Y'] = - Y'$.

Consider the action of $X'$ on the space $V$ of linear forms
 on $\Bbb C^{4}$. Since the
eigenvalues of $X$ are distinct so are the eigenvalues of $X'$. Thus
$V$ decomposes as
$$
V = \bigoplus_{i=0}^{3} V_{\lambda_i} \, ,
$$
where $V_{\lambda_i}$ is the $\lambda_i$-eigenspace for the action
of $X'$, i.e., $X'(v)= \lambda_i v$ for every $v \in V_{\lambda_i}$.
If $v \in V_{\lambda_i}$ then
$$
 [X',Y'](v) = X'(Y'(v)) - Y'(\lambda_i v) = -Y'(v) \, .
$$
and consequently $X'(Y'(v)) = ( \lambda_i -1 ) Y'(v)$, i.e., when
non zero  $Y'(v)$ is an $(\lambda_i -1)$-eigenvector of $X'$.

By semicontinuity,  the rank of $Y'$ is at least $3$ and on the
other hand the equation above  implies that it  is at most $3$.
Moreover we also have that $\ker Y'$ must be equal to one of the
$V_{\lambda_i}$, say $V_{\lambda_k}$. After  reordering, we obtain
 that $\lambda_i= i$.

At this point we have shown that $X'$ is conjugate to $X$. But
$[X',Y']=-Y'$ implies that
$$Y= \sum_{i=0}^{2} \lambda_i z_{i+1}
\frac{\partial}{\partial z_i} \, , $$ where $\lambda _i \in \Bbb C$.
It is now evident that the Lie algebra generated by $X',Y'$ is
conjugate to the one generated by $X,Y$ by an element in
$\text{GL}(\Bbb C^{4})$.
\newline
\newline
(3.4) The second alternative proof is the following.
Let $\bar \omega_H$ be as in (3.1). By assigning to each
differential form its singular set,
the orbit of $\bar \omega_H$ under the automorphism
group of $\Bbb P^3$ maps onto the space of pointed
twisted cubics in $\Bbb P^3$ and hence has dimension at least equal to
$\text{dim Aut}(\Bbb P^3)-\text{dim Aut}(\Bbb P^1)+ 1 = 15-3+1 =
13$.

On the other hand, the tangent space to $\Cal F(3, 2)$ at the point
$\bar \omega_H$ is given  by the forms
$\eta \in H^0(\Bbb P^3, \Omega^1_{\Bbb P^3}(4))$
such that $(\bar \omega_H + \epsilon \eta) \wedge d(\bar \omega_H + \epsilon \eta) = 0$
(modulo $\epsilon^2$), that is:
$$\bar \omega_H  \wedge d \eta + \eta  \wedge d(\bar \omega_H) = 0 $$
($\eta$ defined up to constant multiple of $\bar \omega_H$).
One checks, by hand or by computer,
that the space of solutions $\eta$ of this system of linear
equations has dimension 13. It follows that $\Cal E$ is an
irreducible component of $\Cal F(3, 2)$. Furthermore, the
integrability condition provides $\Cal F(3, 2)$ with a natural
structure of scheme and the tangent space calculation above also
implies that $\Cal E$ is a {\it reduced} component.
\newline
\newline
(3.5) The component $\Cal E$ considered in this article admits some
direct generalizations; let us make some remarks about them.
\newline
\newline
a) For $r \ge 5$, the natural action of the group $PGL(2, \Bbb C)$
on the projective space $\Bbb P^r = \Bbb P S^r(\Bbb C^2)$ of binary
forms of degree $r$ induces a rigid foliation of dimension three and
hence provides an irreducible component of the space $\Cal
F_{r-3}(r, 3)$ of foliations of codimension $r-3$ and degree three
in $\Bbb P^r$ (\cite{CP}, Example (6.6)).
\newline
\newline
Notice that the foliation induced by the action of $PGL(2, \Bbb C)$
on binary forms of degree $r=4$, considered in this article, is not rigid.
This follows from a general fact proved in  \cite{CP}, Proposition (6.5),
or may be seen directly as follows: in (2.6) we observed that the
1-form $\omega$ defining this foliation belongs to the component
$\Cal R(2, 3) \subset \Cal F(4, 3)$ and in fact it is clear that the closure
of the orbit of $\omega$ is a proper subvariety of $\Cal R(2, 3)$
since they have different dimension.
\newline
\newline
b) Let $\text{Aff}(\Bbb C) \subset PGL(2, \Bbb C)$ be the affine
group in one variable. The action of $\text{Aff}(\Bbb C)$ on $\Bbb P
S^r(\Bbb C^2)$ obtained by restricion of the action of $PGL(2, \Bbb
C)$ considered in a) defines a foliation of dimension two. For $r
\ge 4$ these are rigid and define irreducible components of $\Cal
F_{r-2}(r, 2)$ (\cite{CP}, Example (6.8)).
\newline
\newline
c) The component $\Cal E \subset \Cal F(3, 2)$ is the first member
of an infinite family of rigid components $\Cal E(n) \subset \Cal
F(n, n-1)$ defined for $n \ge 3$ in Theorem 4 of \cite{CP}. 
\newline
\newline
It would
be interesting to carry out an analysis of these components similar
to what was done here with the exceptional component $\Cal E$.

\enddocument